\documentclass[11pt]{article}
\usepackage[T1]{fontenc}
\usepackage{lmodern}
\usepackage[a4paper,margin=27mm]{geometry}
\usepackage{amsmath,amssymb,amsthm,microtype}
\usepackage{tikz}
\usepackage{hyperref}
\newtheorem{theorem}{Theorem}
\newtheorem{question}{Question}
\newtheorem{lemma}[theorem]{Lemma}
\title{Edge cuts of bounded degree in sparse graphs}
\author{Anna Neopryatnaya, Vladimir Verdenko, Vsevolod Voronov\footnote{Caucasus Mathematical Center of Adyghe State University; Moscow Institute of Physics and Technology; \\email: v-vor@yandex.ru}}
\date{19 November 2026}

\begin{document}
\maketitle
\begin{abstract}
Every simple graph on \(n\ge2\) vertices with at most \(2n-2\) edges
can be split into two nonempty parts with at most two crossing edges
at each vertex.  For \(d\ge3\)
we bound the minimum number of edges in an \(n\)-vertex graph without a \(d\)-cut.
\end{abstract}

\section{Introduction}

Let \(G\) be a finite simple undirected graph with \(n\ge2\) vertices
and \(m\) edges. For an integer \(d\ge1\), a \emph{\(d\)-cut} is a partition
\(V(G)=A\sqcup B\) into two nonempty sets such that every vertex has at most
\(d\) neighbours in the other set. Equivalently, deleting the edges between
\(A\) and \(B\) removes at most \(d\) edges incident to each vertex. We call
\(A\) and \(B\) the \emph{sides} of the cut.

A disconnected graph has a \(d\)-cut for every \(d\ge1\). More generally,
any cut with at most \(d\) crossing edges is a \(d\)-cut, so edge connectivity
at most \(d\) suffices. In particular, a vertex of degree at most \(d\) gives a singleton side.

A graph with no $d$-cut is $d$-immune~\cite{Feghali}. Let $\mu_d(n)$ be the minimum number of edges of a $d$-immune graph on $n$ vertices. A \(d\)-immune graph must have at least \(2d+1\) vertices, since otherwise both sides can have size at most \(d\). For every \(n\ge2d+1\), the complete graph \(K_n\) is \(d\)-immune because one side of every cut has at least \(d+1\) vertices.

Bonsma, Farley and Proskurowski proved the sharp bound
\(\mu_1(n)=\lceil3(n-1)/2\rceil\) for \(n\ge3\) and
characterized the \(1\)-immune graphs attaining this bound~\cite{BFP}.
For \(d\ge3\), Gomes and Sau proved that every graph with maximum degree $\Delta \leq d+2$ has a $d$-cut, and that for fixed $d\geq 1$ deciding whether a
\(d\)-cut exists is NP-hard even on \((2d+2)\)-regular graphs~\cite{GS}. They also noted that $K_{d+1,2d+1}$ is $d$-immune.


Chernyshev, Rauch and Rautenbach~\cite{CRR} rediscovered the conjecture by Atsushi Kaneko, proposed at 7th C5 Graph Theory Workshop, 2003, that every connected graph with \(m<3n-6\) has a \emph{forest cut}, i.e. a vertex set \(S\) such that \(G-S\) is disconnected and \(G[S]\) is a forest. They proved the planar
case. Botler et al.~\cite{BotlerEtAl} proved the conclusion for
\(m<(9n-15)/4\), and Bogdanov et al.~\cite{BogdanovEtAl} extended it to
\(m<(19n-28)/8\) for \(n\ge4\). The edge \(2\)-cut bound arose from an
attempt to prove a special case of this conjecture.

\begin{theorem}\label{thm:main}
Every graph with \(m\le2n-2\) has a \(2\)-cut. The bound is sharp
for \(n=5,6,8\). For every \(k\ge 1\) there is a 2-immune graph with $7k+1$ vertices and $15k$ edges.
\end{theorem}

\begin{theorem}\label{thm:dhigh}
For every integer \(d\ge3\) and every
\(\varepsilon>0\) there is an \(N=N(d,\varepsilon)\) such that
\[
 \begin{alignedat}{2}
 \mu_d(n)&\ge\left\lfloor\frac{(d+1)(d+3)n}{2(d+2)}\right\rfloor+1
       &\qquad& (n\ge2d+1),\\[2pt]
 \mu_d(n)&\le\left(\frac{d^2+3d+1}{2d+1}+\varepsilon\right)(n-1)
       && (n\ge N).
 \end{alignedat}
\]
\end{theorem}

\section{The \texorpdfstring{\(2\)}{2}-cut threshold}

An \emph{even factor} of \(G\) is a spanning subgraph \(H\) with positive even
degree at every vertex. Its edge complement \(F=G-E(H)\) has
\(d_F(v)\equiv d_G(v)\pmod 2\) and \(d_F(v)\le d_G(v)-2\).
We use the following existence theorem.
\begin{theorem}[Fleischner~\cite{Fleischner}]\label{thm:factor}
Every connected bridgeless graph of minimum degree at least three has
an even factor.
\end{theorem}

Since \(|E(H)|\ge n\), the bound \(m\le2n-2\) gives \(|E(F)|\le n-2\), so \(F\) is disconnected. If \(H\) is \(2\)-regular, each component of \(F\) already gives a \(2\)-cut because all crossing edges lie in \(H\). In general, we merge two parts whenever a vertex in one has at least three neighbours in the other. 


\begin{lemma}\label{lem:merge}
For a partition \(\mathcal P\) of \(V(G)\), let \(b(v)\) count the neighbours
of \(v\) outside its part. Set \(\varphi(\mathcal P)=|\mathcal P|-\sum_v\lfloor b(v)/3\rfloor\).
If \(\varphi(\mathcal P)\ge2\), then \(G\) has a \(2\)-cut.

\end{lemma}
\begin{proof} 
Whenever \(v \in P_1\) has at least three neighbours in $P_2 \neq P_1$, merge
$P_1$ and $P_2$. Denote the resulting partition by $\mathcal P'$. Then $|\mathcal P'|= |\mathcal P|-1$ and \(b(v)\)
decreases by at least three. The merge cannot increase \(b(w)\) for any other vertex \(w\), so
\[
\varphi(\mathcal P') \geq \varphi(\mathcal P) \geq 2.
\]

At termination every vertex outside a given part has at most two neighbours
in it. At most \(\sum_v\lfloor b(v)/3\rfloor\) parts contain a vertex with
\(b(v)\ge3\). Hence at least two parts \(P_1^*,P_2^*\) remain whose vertices
each have at most two outside neighbours. Either part and its complement form
a \(2\)-cut. If only two parts remain, they give the same cut.
\end{proof}

When \(\varphi(\mathcal P)\ge2\), at least two initial parts already satisfy \(b(v)\le2\) at
every vertex they contain. A vertex outside such a part may still have three or more neighbours in it.

\begin{proof}[Proof of the cut threshold in Theorem~\ref{thm:main}]
Assume \(m\le2n-2\). If the edge connectivity is at most two, a cut with
at most two crossing edges suffices. Otherwise \(G\) is bridgeless with minimum degree
at least three, so Theorem~\ref{thm:factor} gives an even factor \(H\). Take the components of
\(G-E(H)\) as the parts of \(\mathcal P\). Since a connected graph on
\(s\) vertices has at least \(s-1\) edges,
\(|\mathcal P|\ge n-m+|E(H)|\ge2\), using \(|E(H)|\ge n\).
Every edge between parts lies in
\(H\). If \(\deg_H(v)=2r\), then \(r\ge1\) and
\(\lfloor b(v)/3\rfloor\le\lfloor2r/3\rfloor\le r-1\). Summing gives
\(\sum_v\lfloor b(v)/3\rfloor\le |E(H)|-n\). Thus
\(\varphi(\mathcal P)\ge2n-m\ge2\), and Lemma~\ref{lem:merge} gives a \(2\)-cut.

For sharpness, \(K_5-e\) has minimum degree three, ruling out a singleton side. In a \(2+3\) partition at least one of the two vertices
on the smaller side is adjacent to all three across. For \(n=6\), consider
the complete tripartite graph \(K_{1,2,3}\). The side opposite its universal
vertex has at most two vertices in any \(2\)-cut, so each would have degree
at most three. The degree-three vertices form an independent set, so each
would have all three neighbours across. In \(K_{3,5}\), at least three
of the five vertices in the larger class lie on one side. In a \(2\)-cut
they force the other class and then all remaining vertices onto that
side~\cite{GS}. Each of the three graphs has \(2n-1\) edges.
\end{proof}

\begin{proof}[Proof of the upper bound in Theorem~\ref{thm:main}]
A \emph{one-vertex sum} identifies one vertex from each of several disjoint
\(2\)-immune graphs. A \(2\)-cut splitting a summand would restrict to a
\(2\)-cut of that summand, so all summands must lie on the side of the common
vertex. Taking \(k\) copies of \(K_{3,5}\) gives \(7k+1\) vertices and \(15k\) edges.

For the \emph{hub construction}, double each edge of a \(2k\)-cycle (\(k\ge2\))
and add a matching between opposite vertices. Subdivide every edge and join \(5k\)
subdivision vertices to a hub \(h\). The resulting simple graph has \(7k+1\) vertices
and \(15k\) edges. 



\end{proof}

\section{The lower bound for \texorpdfstring{\(d\ge3\)}{d >= 3}}

\begin{lemma}\label{lem:local}
If \(d\ge3\) and \(G\) has no \(d\)-cut, then \(\delta(G)\ge d+1\).
The vertices of degree \(d+1\) form an independent set \(I\), and each
vertex of degree \(d+2\) has at most one neighbour in \(I\).
\end{lemma}
\begin{proof}
A set of at most \(d\) vertices meets the cut condition outside. A
vertex of degree at most \(d\) or an adjacent pair of vertices of degree
\(d+1\) also meets it inside. So does a vertex of degree \(d+2\)
together with two neighbours in \(I\), since the centre has two
internal neighbours and each of the others has one.
\end{proof}

\begin{proof}[Proof of the lower bound in Theorem~\ref{thm:dhigh}]
Suppose \(G\) has no \(d\)-cut. Give each vertex \(v\) weight
\[w_0(v) = (d+2)\deg_G(v).\] 

Each vertex outside \(I\) transfers one unit to each of its neighbours in \(I\).
Let \(w_1(v)\) be the resulting weight of \(v\). Put \[K=(d+1)(d+3)=(d+2)^2-1.\]
Every vertex \(v\in I\) receives \(d+1\) units and ends with weight \(w_1(v)=K\). By Lemma~\ref{lem:local},
a vertex of degree \(d+2\) gives away at most one unit and keeps at
least \(K\). A vertex of degree \(k\ge d+3\) keeps at least
\((d+1)k\ge K\). Summing the weights gives \[\sum_{v\in V(G)}w_1(v)=2(d+2)m\ge Kn.\]

Equality would require each vertex outside \(I\) to have degree \(d+2\)
and exactly one neighbour in \(I\), or degree \(d+3\) and all neighbours
in \(I\). If a vertex \(h\) of the second type exists, take \(h\) and
three of its neighbours \(v_1,v_2,v_3\) and put \(A_1=\{h,v_1,v_2,v_3\}\). Each chosen vertex has at most \(d\)
outside neighbours. An outside neighbour of \(h\) belongs to \(I\) and
is adjacent to none of \(v_1,v_2,v_3\), since \(I\) is independent. Every other outside vertex
has at most three neighbours in \(A_1\).

Otherwise \(V(G)\setminus I\) induces a nonempty \((d+1)\)-regular graph \(H\).
Take adjacent vertices \(u_1,u_2\) of \(H\) and their respective neighbours
\(v_1,v_2\in I\), which may coincide. Let \(A_2=\{u_1,u_2,v_1,v_2\}\). Each chosen vertex
has at most \(d\) neighbours outside. An outside vertex has at most one
neighbour in \(I\). Its only other possible neighbours in \(A_2\) are \(u_1,u_2\),
giving at most \(3\le d\) in all.

Both candidate sides \(A_1,A_2\)
are proper since they have at most four vertices and \(n\ge2d+1\).
Either gives a \(d\)-cut. Hence \(2(d+2)m>Kn\), proving the integer bound.
\end{proof}

\section{The upper bound for \texorpdfstring{\(d\ge3\)}{d >= 3}}

For \(d\ge3\), the following construction improves the asymptotic upper bound obtained from
copies of \(K_{d+1,2d+1}\) sharing one vertex. We start with a biregular bipartite
graph of large girth, replace each vertex on one side by a path, and add a vertex
adjacent to every vertex on the other side.

\begin{lemma}[F\"uredi et al.~{\cite[Theorem~A]{FLSUW}}]\label{lem:base}
For integers \(a,b\ge2\), there are finite simple bipartite graphs of
arbitrarily large girth with degrees \(a\) on one side and \(b\) on the
other.
\end{lemma}

\begin{proof}[Proof of the upper bound in Theorem~\ref{thm:dhigh}]
Fix an integer \(t\ge1\). By Lemma~\ref{lem:base}, take a bipartite graph \(Q\)
of girth greater than \(2d\) with degrees \(dt+2\) on \(X\) and \(d+1\)
on \(Y\). Replace each vertex of \(X\) by a path \(P\) on \(t\) vertices.
Assign its \(dt+2\) incident edges to the path vertices, giving each
\(w\) exactly \(d+2-\deg_P(w)\) of them (Figure~\ref{fig:path}).
The path degrees sum to \(2(t-1)\), so all edges are assigned.
Add a vertex \(u\) adjacent to all of \(Y\). In the resulting graph
\(G\), every vertex except \(u\) has degree \(d+2\).

\begin{figure}[htb]
\centering
\begin{tikzpicture}[x=1cm,y=1cm]
\draw[thick] (0,0)--(4,0);
\foreach \x in {0,4}
  \foreach \s in {-.65,-.22,.22,.65}
    \draw[blue!60!black] (\x,0)--({\x+\s},.65) node[circle,fill,inner sep=1pt] {};
\foreach \s in {-.45,0,.45}
  \draw[blue!60!black] (2,0)--({2+\s},.65) node[circle,fill,inner sep=1pt] {};
\foreach \x in {0,2,4}
  \filldraw[fill=white] (\x,0) circle (2.5pt);
\end{tikzpicture}
    \caption{For \(d=t=3\), each path has three vertices and two edges. The endpoints have four neighbours in \(Y\) (filled dots) and the middle vertex has three. All three path vertices have degree five.}
\label{fig:path}
\end{figure}
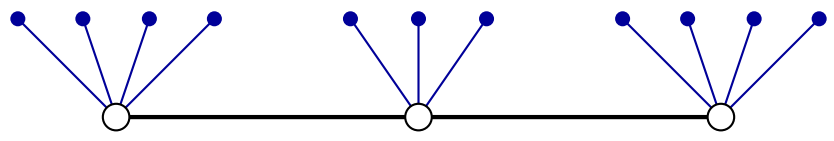

In a \(d\)-cut, the condition at \(u\) gives \(|B\cap Y|\le d\) for the
opposite side \(B\). Each vertex of \(B\) has at least two neighbours
in \(B\). The graph \(Q[X\cup(B\cap Y)]\)
is a forest, since any cycle there would alternate through at most \(d\)
vertices of \(Y\) and have length at most \(2d\). Replacing its
\(X\)-vertices by paths leaves each original edge a bridge. The resulting
forest contains the nonempty graph \(G[B]\), contradicting its minimum degree.

Write \(x=|X|\), \(y=|Y|\), and let the constructed graph have \(n_0\)
vertices and \(m_0\) edges. Counting incidences and degrees gives
\[
 (d+1)y=(dt+2)x,\qquad n_0-1=tx+y,\qquad 2m_0=(d+2)(n_0-1)+y.
\]
Consequently \(m_0=U_d(n_0-1)+x/(2d+1)\), where
\(U_d=(d^2+3d+1)/(2d+1)\). Given \(\varepsilon>0\), choose \(t\) so that
\(1/((2d+1)t)\le\varepsilon/2\). Since \(tx\le n_0-1\), this gives
\(m_0\le (U_d+\varepsilon/2)(n_0-1)\).

For any \(n\ge n_0\), write \(n-1=q(n_0-1)+r\) with \(0\le r<n_0-1\).
Identify the vertices \(u\) of \(q\) disjoint copies of the construction.
Add \(r\) new vertices, joining each to the same \(d+1\) distinct vertices
of the first copy. A \(d\)-cut splitting a copy would restrict to a
\(d\)-cut of that copy. Thus all copies would lie on the side of the
common vertex, forcing each new vertex there through its \(d+1\) neighbours.
The resulting \(n\)-vertex graph has \(qm_0+(d+1)r\) edges, at most
\((U_d+\varepsilon/2)(n-1)+(d+1)(n_0-2)\). Since \(n_0\) is fixed, this
is at most \((U_d+\varepsilon)(n-1)\) for all sufficiently large \(n\).
\end{proof}
\section{Conclusion}


The matching-cut bound is sharp for every \(n\ge3\). For \(2\)-cuts, the bound
appears to be sharp only at small orders. The exhaustive searches give
\(\mu_2(n)=2n\) for \(9\le n\le14\).
For \(n=7k+1\) with \(k\ge2\), both constructions in the proof of
Theorem~\ref{thm:main} give \(m=15(n-1)/7\) edges. 

\begin{question}
Is \(\mu_2(n) \geq 15(n-1)/7\) for $n\geq 5$?
\end{question}

For \(d \geq 3\), it would also be interesting to determine the limiting density $\alpha_d = \lim_{n \to \infty}\frac{\mu_d(n)}{n-1}$ at which \(d\)-immune graphs exist.

\paragraph{AI use statement.} GPT-6 Astra was used in searching for proofs and in editing the paper. The proof for 2-cuts was found in the course of lengthy discussions in which some of the ideas belonged to the authors and the other part to AI. The proofs were verified by the authors, and some arguments were significantly simplified.

\end{document}